\documentclass{amsart}
\usepackage{amsmath}
\usepackage{epic}
\usepackage{eepic}

\newtheorem{thm}{Theorem}[section]
\newtheorem{lem}[thm]{Lemma}

\newtheorem{clm}[thm]{Claim}

\theoremstyle{definition}
\newtheorem{defn}[thm]{Definition}

\newcommand{\Z}{{\bf{Z}}}

\newcommand{\cC}{{\mathcal{C}}}
\newcommand{\cG}{{\mathcal{G}}}
\newcommand{\cH}{{\mathcal{H}}}
\newcommand{\cI}{{\mathcal{I}}}

\newcommand{\cR}{{\mathcal{R}}}

\newcommand{\co}{\colon\thinspace}

\begin{document}

\title{Braid groups are linear}                    
\author{Stephen J. Bigelow}                  
\address{University of California, Berkeley, California 94720} 
\email{bigelow@math.berkeley.edu}


\begin{abstract}
The braid groups $B_n$ can be defined as
the mapping class group of the $n$-punctured disc.
The Lawrence-Krammer representation of the braid group
is the induced action on
a certain twisted second homology
of the space $C$ of unordered pairs of points
in the $n$-punctured disc.
Recently, Daan Krammer showed that this is a faithful
representation in the case $n = 4$.
In this paper, we show that it is faithful for all $n$.
\end{abstract}

\maketitle

\section{Introduction}
\label{INTRODUCTION}

Let $B_n$ denote Artin's braid group on $n$ strands.
Recently, Krammer \cite{krammer}
proved that a certain representation of the braid groups
is faithful on $B_4$.
The representation he used is essentially the same
as one used by Lawrence in \cite{lawrence}
to give a topological definition of 
a certain summand of the Jones representation.
We call this representation
the Lawrence-Krammer representation.
In this paper, we prove the following.

\begin{thm}
The Lawrence-Krammer representation of $B_n$
is faithful for all $n$.
\end{thm}

This proves that braid groups are linear,
thus solving a long-standing open problem.
Our proof can be seen as a sort of converse
to the construction of elements of the kernel of the Burau representation
given in \cite{moody}, \cite{long-paton} and \cite{bigelow}

\subsection{Definitions}

Let $D$ be an oriented disc in the complex plane.
Fix a set $P \subset D$ 
consisting of $n$ distinct points $p_1,\dots,p_n$
in the interior of $D$.
Let $\cH(D,P)$ be the group of all
homeomorphisms $h \co D \rightarrow D$
such that $h(P) = P$ and $h$ fixes $\partial D$ pointwise.
Let $\cI(D,P)$ be the group of all such homeomorphisms which are
isotopic to the identity relative to $P \cup \partial D$.
We define the braid group $B_n$ to be the group $\cH(D,P)/\cI(D,P)$.
See \cite{birman} for other equivalent definitions of these groups
and a good introduction to their basic properties.

Let $C$ denote the space of all unordered pairs of distinct points
in $D \setminus P$.
In other words,
\[
C = 
\frac
{((D \setminus P) \times (D \setminus P)) \setminus \{(x,x)\}}
{(x,y) \sim (y,x)}
.
\]
Let $d_1$ and $d_2$ be distinct points in $\partial D$.
Let $c_0 = \{d_1,d_2\}$ be a basepoint for $C$.

We now define a map $\phi$
from $\pi_1(C,c_0)$ to the free Abelian group
with basis $\{q,t\}$.
Let $\alpha$ be a closed curve in $C$ based at $c_0$
representing an element $[\alpha]$ of $\pi_1(C,c_0)$.
We can write $\alpha$ in the form
\[
\alpha(s) = \{\alpha_1(s),\alpha_2(s)\}
\]
for some arcs $\alpha_1$ and $\alpha_2$ in $D \setminus P$.
Let
\[
a = 
\frac{1}{2\pi i}
\sum_{j=1}^n
\left(
\int_{\alpha_1} \frac{dz}{z-p_j} + \int_{\alpha_2} \frac{dz}{z-p_j}
\right).
\]
Let 
\[
b = \frac{1}{\pi i} \int_{\alpha_1 - \alpha_2} \frac{dz}{z}.
\]
Let $\phi([\alpha]) = q^a t^b$.

This definition requires some explanation.
If $\alpha_1$ and $\alpha_2$ are closed loops then
$a$ is the sum of the winding numbers of $\alpha_1$ and $\alpha_2$
around each of the puncture points,
and $b$ is twice the winding number of $\alpha_1$ and $\alpha_2$
around each other.
However $\alpha_1$ and $\alpha_2$ are not necessarily closed loops,
but may ``switch places''.
In this case,
$a$ is the sum of the winding numbers
of the closed loop $\alpha_1\alpha_2$
around each of the puncture points,
and $\alpha_1 - \alpha_2$ satisfies
$$(\alpha_1 - \alpha_2)(1) = - (\alpha_1 - \alpha_2)(0),$$
which implies that $b$ is an odd integer.

Let $\tilde{C}$ be the covering space of $C$
whose fundamental group is the kernel of $\phi$.
Choose a lift $\tilde{c}_0$ of $c_0$ to $\tilde{C}$.
Let $\Lambda$ denote the ring $\Z[q^{\pm 1},t^{\pm 1}]$.
The homology group $H_2(\tilde{C})$
can be considered as a $\Lambda$-module,
where $q$ and $t$ act by covering transformations.

Let $\sigma \in B_n$.
There is an induced action of $\sigma$ on $C$.
Let $\tilde{\sigma}$
be the lift of $\sigma$ to a map from $\tilde{C}$ to itself
which fixes $\tilde{c}_0$.
This induces a $\Lambda$-module automorphism $\tilde{\sigma}_*$
of $H_2(\tilde{C})$.
The Lawrence-Krammer representation is the map
from $B_n$ to $GL(H_2(\tilde{C}))$
taking $[\sigma]$ to $\tilde{\sigma}_*$.

\subsection{Outline}

In Section \ref{fn}
we will define forks and noodles.
These will be one-dimensional objects in the disc
designed to represent
elements of the second homology and cohomology of $\tilde{C}$.
We define a pairing between forks and noodles,
which will be preserved by any element of the kernel
of the Lawrence-Krammer representation.

In Section \ref{key}
we prove that pairing between forks and noodles
detects geometric intersection
between the corresponding edges in the disc.
We use this to show that a braid in the kernel of
the Lawrence-Krammer representation
must be trivial.

In Section \ref{lkmatrices}
we compute the Lawrence-Krammer representation explicitly 
in terms of generators and basis elements.

\subsection{Notation}

If $\alpha$ and $\beta$ are arcs in $D \setminus P$
such $\alpha(s) \neq \beta(s)$ for all $s \in I$
then we define $\{\alpha,\beta\}$ to be the arc in $C$
given by
$$\{\alpha,\beta\}(s) = \{\alpha(s),\beta(s)\}.$$
If $y$ is a point in $D \setminus P$
and $\alpha$ is an arc in $D \setminus (P \cup \{y\})$
then we define $\{\alpha,y\}$ to be the arc in $C$
given by
$$\{\alpha,y\}(s) = \{\alpha(s),y\}.$$
The same arc can be denoted by $\{y,\alpha\}$.

If $x$ and $y$ are elements of a group then we use the notation
$$x^y = y^{-1}xy$$ 
and 
$$[x,y] = x^{-1}y^{-1}xy.$$
Throughout this paper, $I$ will denote the interval $[0,1]$.
Braids compose from right to left.
Arcs compose from left to right.
\section{Forks and Noodles}
\label{fn}

In this section
we define forks and noodles and a pairing between them.
The idea of using a fork
to represent an element of $H_2(\tilde{C})$
is due to Krammer \cite{krammer}.

A {\em fork} is an embedded tree $F \subset D$
with four vertices $d_1$, $p_i$, $p_j$ and $z$
such that
$F \cap \partial D = \{d_1\}$,
$F \cap P = \{p_i,p_j\}$, and
all three edges have $z$ as a vertex.
The edge containing $d_1$ is called the {\em handle} of $F$.
The union of the other two edges is a single edge,
which we call the {\em tine edge} of $F$ and denote by $T(F)$.

\begin{figure}
\centering

\begin{picture}(120,130)(0,-10)

\thinlines
\put(60,60){\ellipse{120}{120}}

\drawline(50,1)(50,70)
\qbezier(20,60)(30,70)(40,70)
\drawline(40,70)(80,70)
\qbezier(80,70)(90,70)(100,60)
\put(56,74){$T(F)$}
\put(38,20){$F$}

\drawline(70,1)(70,50)
\qbezier(20,60)(30,50)(40,50)
\drawline(40,50)(80,50)
\qbezier(80,50)(90,50)(100,60)
\put(56,54){$T(F')$}
\put(70,20){$F'$}

\put(45,-10){$d_1$}
\put(65,-10){$d_2$}

\put(20,60){\circle*{4}}
\put(100,60){\circle*{4}}

\put(49,90){\circle*{4}}
\put(60,90){\circle*{4}}
\put(71,90){\circle*{4}}
\put(54.5,100){\circle*{4}}
\put(65.5,100){\circle*{4}}
\end{picture}
\caption{A fork $F$ and a parallel copy $F'$.}
\label{doublefork}
\end{figure}
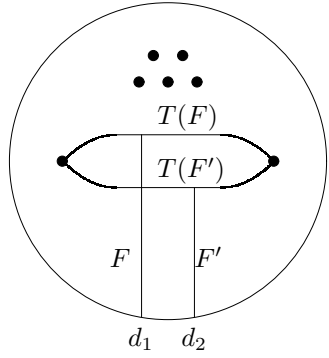
For a given fork $F$
we can define a {\em parallel copy} of $F$
to be an embedded tree $F'$
as shown in Figure \ref{doublefork}.
The five punctures at the top of the figure
may be replaced by any number,
and any orientation-preserving self-homeomorphism
may be applied to the entire figure.
The tine edge $T(F')$ of $F'$ is defined analogously to that of $F$.

For any fork $F$,
we define a surface $\tilde{\Sigma}(F)$ in $\tilde{C}$ as follows.
Let $F'$ be a parallel copy of $F$.
Let $z$ be the vertex contained in all three edges of $F$,
and let $z'$ be the vertex contained in all three edges of $F'$.
We define the surface $\Sigma(F)$ 
to be the set of points in $C$
which can be written in the form $\{x,y\}$,
where $x \in T(F) \setminus P$ 
and $y \in T(F') \setminus P$.
Let $\beta_1$ be an arc from $d_1$ to $z$
along the handle of $F$
and let $\beta_2$ be an arc from $d_2$ to $z'$
along the handle of $F'$.
Let $\tilde{\beta}$ be the lift of 
$\{\beta_1,\beta_2\}$ to $\tilde{C}$
which starts at $\tilde{c}_0$.
Let $\tilde{\Sigma}(F)$ be the lift of $\Sigma(F)$ to $\tilde{C}$
which contains $\tilde{\beta}(1)$.

A {\em noodle}
is an embedded edge $N \subset D \setminus P$
with endpoints $d_1$ and $d_2$.
For a given noodle $N$
we define the surface $\Sigma(N)$ to be the set of points $\{x,y\} \in C$
such that $x$ and $y$ are distinct points in $N$.
Let $\tilde{\Sigma}(N)$ be the lift of $\Sigma(N)$ to $\tilde{C}$
which contains $\tilde{c}_0$.

Let $F$ be a fork and let $N$ be a noodle.
We define $\langle N,F \rangle \in \Lambda$ as follows.
\[
\langle N,F \rangle = 
\sum_{a,b \in Z}q^a t^b 
(q^a t^b \tilde{\Sigma}(N), \tilde{\Sigma}(F)).
\]
Here $q^a t^b \tilde{\Sigma}(N)$ denotes
the image of $\tilde{\Sigma}(N)$
under the covering transformation $q^a t^b$,
and $(q^a t^b \tilde{\Sigma}(N), \tilde{\Sigma}(F))$
denotes the algebraic intersection number
of the two surfaces in $\tilde{C}$.

\begin{lem}[The Basic Lemma]
\label{blem}
The pairing between forks and noodles is well-defined.
Furthermore, if $[\sigma]$ lies in 
the kernel of the Lawrence-Krammer representation then 
$$\langle N,F \rangle = \langle N,\sigma(F) \rangle$$
for every fork $F$ and noodle $N$.
\end{lem}

\subsection{Proof of the Basic Lemma}
\label{BLEM}

The problem is that
one cannot necessarily define an algebraic intersection number
between two properly embedded surfaces,
since it might be possible
to eliminate intersections
by pushing them off to infinity.
We overcome this problem by proving
the existence of an immersed closed surface $\tilde{\Sigma}_2(F)$
which is equal to $(1-q)^2(1+qt)\tilde{\Sigma}(F)$
outside a small neighborhood of
the punctures.

Let $F$ be a fork.
Let the endpoints of $T(F)$ be $p_i$ and $p_j$.
Let $\nu(p_i)$ and $\nu(p_j)$ be disjoint
$\epsilon$-neighborhoods of $p_i$ and $p_j$ respectively
such that $\nu(p_k) \cap P = \{p_k\}$ for $k=i,j$.
Let $U$ be the set of $\{x,y\} \in C$
such that at least one of $x$ and $y$
lies in $\nu(p_i) \cup \nu(p_j)$.
Fix a basepoint $u_0 = \{u_1,u_2\} \in U$,
where $u_1 \in \nu(p_i)$ and $u_2 \in \nu(p_j)$.
Let $\tilde{U}$ be the pre-image of $U$ in $\tilde{C}$.
Choose a lift $\tilde{u}_0$ of $u_0$ to $\tilde{C}$.
We start by analyzing $\pi_1(U,u_0)$.

\begin{figure}
\centering

\begin{picture}(280,130)(4,0)

\thinlines

\put(60,60){\ellipse{120}{120}}
\put(55,125){$\nu(p_i)$}

\put(220,60){\ellipse{120}{120}}
\put(215,125){$\nu(p_j)$}

\drawline(85,60)(105,40)
\put(95,50){\vector(1,-1){0}}
\put(95,50){$\alpha_1$}

\drawline(105,40)(175,40)
\put(140,40){\vector(1,0){0}}
\put(138,33){$\alpha_2$}

\drawline(175,40)(195,60)
\put(185,50){\vector(1,1){0}}
\put(174,52){$\alpha_3$}

\drawline(195,60)(175,80)
\put(185,70){\vector(-1,1){0}}
\put(189,73){$\beta_1$}

\drawline(175,80)(105,80)
\put(140,80){\vector(-1,0){0}}
\put(140,85){$\beta_2$}

\drawline(105,80)(85,60)
\put(95,70){\vector(-1,-1){0}}
\put(86,75){$\beta_3$}

\put(60,60){\circle*{4}}
\put(50,54){$p_i$}

\put(77,50){$u_1$}

\qbezier(40,60)(40,80)(85,60)
\qbezier(40,60)(40,40)(85,60)
\put(40,58){\vector(0,-1){0}}
\put(30,58){$\gamma_1$}

\put(192,50){$u_2$}

\put(220,60){\circle*{4}}
\put(223,56){$p_j$}

\qbezier(195,60)(240,80)(240,60)
\qbezier(195,60)(240,40)(240,60)
\put(240,62){\vector(0,1){0}}
\put(241,58){$\gamma_2$}

\end{picture}
\caption{Arcs.}
\label{u}
\end{figure}
Using the arcs shown in Figure \ref{u},
we define the following elements of $\pi_1(U,u_0)$.
\begin{eqnarray*}
a_1 &=& \{\gamma_1,u_2\},\\
a_2 &=& \{u_1,\gamma_2\},\\
b_1 &=& \{\alpha_1, \beta_1\beta_2\beta_3\}\{\alpha_2\alpha_3,u_1\},\\
b_2 &=& \{\alpha_1\alpha_2\alpha_3, \beta_1\}\{u_2,\beta_2\beta_3\}.
\end{eqnarray*}
Note that $b_1$ and $b_2$ are homotopic in $C$, but not in $U$.

The following relations hold in $\pi_1(U)$.
\begin{eqnarray}
\label{urel1}  {[}a_1,a_2] &=& 1,\\
\label{urel2}  {[}a_1,b_1 a_1 b_1] &=& 1,\\
\label{urel3}  {[}a_2,b_2 a_2 b_2] &=& 1.
\end{eqnarray}
The first of these is obvious.
The second follows from the fact that
$b_1 a_1 b_1$ is equal in $\pi_1(U)$
to $\{u_1,\delta\}$,
where $\delta$ is a curve based at $u_2$
which passes counterclockwise around $p_1$ and $u_1$.
The third follows by a similar argument.

We now analyze $\pi_1(\tilde{U},\tilde{u}_0)$.
Let $i \co U \to C$ be the inclusion map
and let $i_*$ be the induced map on fundamental groups.
Then $\pi_1(\tilde{U})$ is the kernel of the map $\phi i_*$.
We define the following elements of $\pi_1(\tilde{U},\tilde{u}_0)$.
\begin{eqnarray*}
a &=& a_2^{-1}a_1,\\
b &=& b_2^{-1}b_1,\\
c &=& a_1^{-1}b_1^{-1}a_1b_1,\\
d &=& a_2^{-1}b_2^{-1}a_2b_2.
\end{eqnarray*}
If $x \in \pi_1(\tilde{U},\tilde{u}_0)$
and $y \in \pi_1(U,u_0)$
then the conjugate $x^y = y^{-1}xy$
is also an element of $\pi_1(\tilde{U},\tilde{u}_0)$.
The following relations hold in $\pi_1(\tilde{U},\tilde{u}_0)$.
\begin{eqnarray}
\label{Urel1}  a^{a_1} &=& a,\\
\label{Urel2}  c^{b_1 a_1} c &=& 1,\\
\label{Urel3}  d^{b_2 a_2} d &=& 1,\\
\label{Urel4}  dba^{b_1} &=& ab^{a_1}c.
\end{eqnarray}
To see this,
rewrite these relations in terms of
$a_1$, $a_2$, $b_1$ and $b_2$.
The first three translate into
equations (\ref{urel1}) to (\ref{urel3}).
The fourth translates into a trivial identity.

If $x \in \pi_1(\tilde{U})$,
let $[x]$ denote the corresponding element of $H_1(\tilde{U})$.
Note that
if $x \in \pi_1(\tilde{U})$ and $y \in \pi_1(U)$ then
$[x^y] = \phi(y)^{-1}[x]$.
The relations given in equations (\ref{Urel1}) to (\ref{Urel4})
give rise to the following relations in $H_1(\tilde{U})$.
\begin{eqnarray*}
\label{hrel1} (q^{-1}-1)[a] &=& 0,\\
\label{hrel2} (q^{-1}t^{-1}+1)[c] &=& 0,\\
\label{hrel3} (q^{-1}t^{-1}+1)[d] &=& 0,\\
\label{hrel4} (q^{-1}-1)[b] &=& (q^{-1}-1)[a] - [c] + [d].
\end{eqnarray*}
Combining these relations, we obtain the following.
$$(1-q)^2(1+qt) [b] = 0.$$

Let $[\tilde{\Sigma}(F)]$ be the element of $H_2(\tilde{C},\tilde{U})$
represented by $\tilde{\Sigma}(F)$.
The long exact sequence of relative homology
gives us the following exact sequence of $\Lambda$-modules
\[
H_2(\tilde{C}) \stackrel{j_*}{\to}
H_2(\tilde{C},\tilde{U}) \stackrel{\partial}{\to}
H_1(\tilde{U}).
\]
But $\partial[\tilde{\Sigma}(F)] = [b]$.
It follows that 
$$(1-q)^2(1+qt)[\tilde{\Sigma}(F)] = j_*[\tilde{\Sigma}_2(F)]$$
for some $[\tilde{\Sigma}_2(F)] \in H_2(\tilde{C})$
represented by some closed immersed surface $\tilde{\Sigma}_2(F)$.

Let $N$ be a noodle.
Choose $\nu(p_i)$ and $\nu(p_j)$ small enough
so as not to intersect $N$.
Then
\[
(1-q)^2(1+qt) \langle N,F \rangle =
\sum_{a,b \in \Z}q^at^b(q^at^b\tilde{\Sigma}(N),\tilde{\Sigma}_2(F)).
\]
Now $(q^at^b\tilde{\Sigma}(N),\tilde{\Sigma}_2(F))$
is the algebraic intersection number
between a properly embedded surface
and an immersed closed surface,
so is well-defined.

Suppose $\sigma$ is an element of the kernel of
the Lawrence-Krammer representation.
Then $\sigma(\tilde{\Sigma}_2(F))$
and $\tilde{\Sigma}_2(F)$
represent the same element of homology,
so have the same algebraic intersection with
any $q^at^b\tilde{\Sigma}(N)$.
Thus $\langle N,\sigma(F) \rangle = \langle N,F \rangle$.

\subsection{Alternative proofs}

There are many possible approaches to proving the Basic Lemma.
The proof given above is a compromise of sorts,
since it proves the existence of
an appropriate element of $H_2(\tilde{C})$,
but does so in a non-constructive way.

It is possible to give a more constructive proof
which uses an explicit computation of $H_2(\tilde{C})$.
One obtains a concrete description of
an immersed genus two surface
which can be seen to be the same as $(1-q)^2(1+qt)\tilde{\Sigma}(F)$
away from the puncture points.
This is perhaps best done in private,
since the details are only convincing
to the person who figures them out.
Some details of a computation of $H_2(\tilde{C})$ 
will be given in Section \ref{lkmatrices}.
See also \cite{lawrence},
where similar methods are used to calculate
the middle homology of the space of 
{\em ordered} $k$-tuples of distinct points in 
the $n$-times punctured disc,
where $k$ can be any positive integer.

It is tempting to seek a less constructive proof
which makes no reference to $\tilde{\Sigma}_2(F)$.
It is intuitively obvious that
the problem of pushing intersections off to infinity
does not arise in the context of forks and noodles.
However this line of reasoning runs into some technical difficulties 
which I feel distract from the true nature of the problem at hand.
A proof that $B_n$ acts faithfully on $H_2(\tilde{C})$
should refer to an element of $H_2(\tilde{C})$.

It is possible to prove that braid groups are linear
without reference to $C$, let alone $H_2(\tilde{C})$.
The Lawrence-Krammer representation can be defined to be
the action of $B_n$ on a $\Lambda$-module
consisting of formal linear combinations of forks
subject to certain relations,
as described by Krammer in \cite{krammer}.
The pairing $\langle N,F \rangle$
can be defined solely in terms of winding numbers.
One must check that this pairing respects
the relations between forks.
The Basic Lemma then follows immediately.
The rest of the proof that 
the Lawrence-Krammer representation is faithful 
proceeds virtually unchanged.
\section{The representation is faithful}
\label{key}

In this section, we prove that 
the Lawrence-Krammer representation is faithful.
The main ingredient in the proof is the following lemma.

\begin{lem}[The Key Lemma]
Let $N$ be a noodle and let $F$ be a fork.
Then $\langle N,F \rangle = 0$
if and only if
$T(F)$ is isotopic relative to $\partial D \cup P$
to an arc which is disjoint from $N$.
\end{lem}

\subsection{Computing the pairing}

We now describe how to compute the pairing 
between a give noodle and fork.
Let $N$ be a noodle and let $F$ be a fork.
By applying a preliminary isotopy, we can assume that
$T(F)$ intersects $N$ transversely at a finite number of points,
which we label $z_1,\dots,z_l$.
Let $F'$ be a parallel copy of $F$.
Choose $F'$ so that $T(F')$ intersects $N$ transversely
at $z'_1,\dots,z'_l$,
where $z_i$ and $z'_i$ are joined
by a short arc in $N$ which meets no other $z_j$ or $z'_j$.
Let $z$ be the vertex which lies in all three edges of $F$,
and let $z'$ be the vertex which lies in all three edges of $F'$.

For every pair $\{z_i,z'_j\}$
there exist unique integers $a_{i,j}$ and $b_{i,j}$
such that $\tilde{\Sigma}(F)$
intersects $q^{a_{i,j}}t^{b_{i,j}}\tilde{\Sigma}(N)$
at a point in the fiber of $\{z_i,z'_j\}$.
Let $\epsilon_{i,j}$ denote the sign of that intersection.
Let $m_{i,j} = q^{a_{i,j}}t^{b_{i,j}}$.
Then
\begin{equation}
\label{pairsum}
\langle N,F \rangle = \sum_{i=1}^l \sum_{j=1}^l \epsilon_{i,j}m_{i,j}.
\end{equation}

We can compute $m_{i,j}$ as follows.
Define the following embedded arcs in $D \setminus P$.
\begin{itemize}
\item $\alpha_1$ from $d_1$ to $z$ along the handle of $F$,
\item $\alpha_2$ from $d_2$ to $z'$ along the handle of $F'$,
\item $\beta_1$ from $z$ to $z_i$ along $T(F)$,
\item $\beta_2$ from $z'$ to $z'_j$ along $T(F')$,
\item $\gamma_1$ from $z_i$ to $d_k$ along $N$,
      where $k=1,2$ is such that $\gamma_1$ does not pass through $z'_j$,
\item $\gamma_2$ from $z'_j$ to $d_{k'}$ along $N$,
      where $k'=1,2$ is such that $\gamma_2$ does not pass through $z_i$.
\end{itemize}
Let 
\[
\delta_{i,j} = 
\{\alpha_1,\alpha_2\}
\{\beta_1,\beta_2\}
\{\gamma_1,\gamma_2\}.
\]
Then
$$m_{i,j} = \phi(\delta_{i,j}).$$

We can calculate the exponent $a_{i,j}$ even more explicitly.
For $k=1,\dots,l$,
let $\zeta_k$ be the arc from $d_1$ to $z_k$ along $F$,
then back to $d_1$ along $N$.
Let $a_k$ be the sum of the winding numbers
of $\zeta_k$ around each of the puncture points.
Let $\zeta$ be the arc from $d_1$ to $d_2$ along $N$,
and then back to $d_1$ moving clockwise 
along $\partial D$.
Let $a$ be the sum of the winding numbers
of $\zeta$ around the puncture points.

\begin{clm}
\label{a}
$a_{i,j} = a_i + a_j + a$.
\end{clm}

\begin{proof}
Let $\zeta'_k$ be the arc from $d_2$ to $z'_k$ along $F'$,
then back to $d_2$ along $N$.
Let $a'_k$ be the sum of the winding numbers
of $\zeta'_k$ around each of the puncture points.
If $z_i$ is closer to $d_1$ than $z'_j$ along $N$
then $a_{i,j} = a_i + a'_j$.
If not, then
$a_{i,j} = a'_i + a_j$.
But $a'_k = a_k + a$ for all $k=1,\dots,l$,
so in either case $a_{i,j} = a_i + a_j + a$.
\end{proof}

In order to compute $\epsilon_{i,j}$
we need to choose orientations for $\Sigma(F)$, $\Sigma(N)$ and $C$.
We choose the orientation on $C$
induced by the product orientation on $D \times D$.
We orient $\Sigma(F)$ and $\Sigma(N)$ as follows.

Let $f \co T(F) \to I$
and $f' \co T(F') \to I$ be diffeomorphisms
which give $T(F)$ and $T(F')$ parallel orientations.
Let $f_1$ and $f_2$ be the maps from $\Sigma(F)$ to $I$
such that if $x \in T(F) \setminus P$ and $y \in T(F') \setminus P$
then $f_1(\{x,y\}) = f(x)$ and $f_2(\{x,y\}) = f'(y)$.
Then the $2$-form $df_1 \wedge df_2$
defines an orientation on $\Sigma(F)$.

Let $g \co N \to I$ be a diffeomorphism
such that $g(d_1) = 0$ and $g(d_2) = 1$.
Let $g_1$ and $g_2$ be the maps from $\Sigma(N)$ to $I$
given by $g_1(\{x,y\}) = \min(g(x),g(y))$
and $g_2(\{x,y\}) = \max(g(x),g(y))$.
Then the $2$-form $dg_1 \wedge dg_2$
defines an orientation on $\Sigma(N)$.

\begin{clm}
\label{sign}
With these orientations,
$\epsilon_{i,j} = -m_{i,i}m_{i,j}m_{j,j}|(q=1,t=-1)$.
\end{clm}

\begin{proof}
By definition, $\epsilon_{i,j}$
is the sign of the volume form 
$dg_1 \wedge dg_2 \wedge df_1 \wedge df_2$
at the point $\{z_i,z'_j\}$ in $C$.
This is determined by the following three things:
\begin{itemize}
\item the sign of the intersection of $N$ with $T(F)$ at $z_i$,
\item the sign of the intersection of $N$ with $T(F')$ at $z'_j$,
\item which of $z_i$ and $z'_j$ is closer to $d_1$ along $N$.
\end{itemize}
These in turn are determined by the following three values
respectively:
\begin{itemize}
\item $m_{i,i}|(q=1,t=-1)$,
\item $m_{j,j}|(q=1,t=-1)$,
\item $m_{i,j}|(q=1,t=-1)$.
\end{itemize}
To see this, observe that in general
the sign of $m_{i',j'}|(q=1,t=-1)$ is determined by whether
the two points switch places in the path $\delta_{i',j'}$.
This in turn is determined by which of $z_{i'}$ and $z'_{j'}$
is closer to $d_1$ along $N$.

It follows from the above considerations
that $\epsilon_{i,j}$ is determined by the value of the product
$m_{i,i}m_{i,j}m_{j,j}$ evaluated at $q=1$ and $t=-1$.
It remains to check that the sign is as claimed.
This can be done by calculating $\epsilon_{i,j}$
in a specific example.
\end{proof}


\subsection{Proof of the Key Lemma}

Let $N$ be a noodle and let $F$ be a fork.
By applying a preliminary isotopy, we can assume that
$T(F)$ intersects $N$ transversely at a finite number of points,
which we label $z_1,\dots,z_l$.
Further, we can assume that $l$
is the minimal possible number of points of intersection.
Let $F'$ be a parallel copy of $F$.
Choose $F'$ so that
and $T(F')$ intersects $N$ transversely
at $z'_1,\dots,z'_l$,
where $z_i$ and $z'_i$ are joined
by a short arc in $N$ which meets no other $z_j$ or $z'_j$.

If $l = 0$ then clearly $\langle N,F \rangle = 0$.
We now assume $l > 0$ and show that $\langle N,F \rangle \neq 0$.
We use the following lexicographic ordering
on the set of monomials $q^a t^b$.

\begin{defn}
We say $q^a t^b \le q^{a'}t^{b'}$
if and only if either
\begin{itemize}
\item $a < a'$, or
\item $a = a'$ and $b \le b'$.
\end{itemize}
For $i,j \in \{1,\dots,l\}$
we say that $m_{i,j}$ is {\em maximal}
if $m_{i,j} \ge m_{i',j'}$
for all $i',j' \in \{1,\dots,l\}$.
\end{defn}

\begin{clm}
\label{keyclaim}
If $m_{i,j}$ is maximal
then $m_{i,i} = m_{j,j} = m_{i,j}$.
\end{clm}

\begin{proof}
Suppose $m_{i,j}$ is maximal.
Then $a_{i,j}$ is maximal among all the integers $a_{k,l}$.
By Claim \ref{a}
it follows that $a_i$ and $a_j$ are maximal
among all the integers $a_{k}$.
Thus $a_{i,i} = a_{j,j} = a_{i,j}$.

We now show that $b_{i,i} = b_{i,j}$.
Suppose not.
Then $b_{i,i} < b_{i,j}$.
Let $\alpha$ be an embedded arc
from $z'_i$ to $z'_j$ along $T(F')$.
Let $\beta$ be an embedded arc 
from $z'_j$ to $z'_i$ along $N$.

If $\beta$ does not pass through the point $z_i$,
let $\delta = \alpha \beta$
and let $w$ be the winding number of $\delta$ around $z_i$.
Then $b_{i,j} - b_{i,i} = 2w$.

If $\beta$ does pass through $z_i$,
first modify $\beta$ in a small neighborhood of $z_i$
so that $z_i$ lies to its left.
Next let $\delta = \alpha \beta$
and let $w$ be the winding number of $\delta$ around $z_i$.
Then $b_{i,j} - b_{i,i} = 2w-1$.

In either case, $w$ is greater than zero.

Let $D_1 = D \setminus \{z_i\}$.
Let $\tilde{D}_1$ be the universal cover of $D_1$.
Let $\tilde{\alpha}$ be a lift of $\alpha$ to $\tilde{D}_1$.
Let $\tilde{\beta}$ be the lift of $\beta$ to $\tilde{D}_1$
which starts at $\tilde{\alpha}(1)$.
Let $\tilde{\gamma}$ be an embedded arc in $\tilde{D}_1$
from $\tilde{\beta}(1)$ to $\tilde{\alpha}(0)$
which intersects $\tilde{\alpha}$ and $\tilde{\beta}$
only at its endpoints.
Let $\gamma$ be the projection of $\tilde{\gamma}$ to $D_1$.
Choose $\tilde{\gamma}$ so that $\gamma$
does not wind around any puncture points.

Let $\tilde{z}'_k$ be the first point on $\tilde{\alpha}$
which intersects $\tilde{\beta}$
(possibly $\tilde{\alpha}(1)$).
This lies in the fiber over $z'_k$ for some $k=1,\dots,l$.
Let $\tilde{\alpha}'$ be the initial segment of $\tilde{\alpha}$
ending at $\tilde{z}'_k$.
Let $\tilde{\beta}'$ be the final segment of $\tilde{\beta}$
starting at $\tilde{z}'_k$.
Let $\tilde{\delta}' = \tilde{\alpha}' \tilde{\beta}' \tilde{\gamma}$.

Now $\tilde{\delta}'$ is a simple closed curve in $\tilde{D}_1$,
so by the Jordan curve theorem it must bound a disc $\tilde{B}$.
Since $\gamma$ passes clockwise around $z_i$,
there is a non-compact region to the right of $\tilde{\delta}'$,
so $\tilde{\delta}'$ must pass counterclockwise around $\tilde{B}$.

Let $\delta'$ be the projection of $\tilde{\delta}'$ to $D$.
Then $a_k - a_i$ is equal to
the sum of the winding numbers of $\delta'$
around each of the points in $P$.
This is equal to the number of points in $\tilde{B}$
which are lifts of a point on $P$.
This must be greater than zero,
since otherwise we could isotope $T(F')$ 
so as to have fewer points of intersection with $N$.
Thus $a_k - a_i$ is greater than zero,
contradicting the fact that $a_i$ is maximal
among all integers $a_{i'}$.
Therefore our assumption that $b_{i,j} > b_{i,i}$ must have been false,
so $b_{i,j} = b_{i,i}$.

The proof that $b_{i,j} = b_{j,j}$ is similar.
This completes the proof of the claim.
\end{proof}

This claim,
together with the formula for $\epsilon_{i,j}$ given in Claim \ref{sign},
implies that if $m_{i,j}$ is maximal
then $\epsilon_{i,j} = -m_{i,j}|(q=1,t=-1)$.
Thus all maximal monomials 
occur with the same sign in equation (\ref{pairsum}).
Therefore $\langle N,F \rangle$ cannot equal zero.
This completes the proof of the Key Lemma.

\subsection{Proof of the Theorem}

In this subsection we use
the Basic Lemma and the Key Lemma
to prove that the Lawrence-Krammer representation is faithful.
The following lemma will be useful.

\begin{lem}
\label{digon}
Let $\alpha$ and $\beta$ be simple closed curves in $D \setminus P$
which intersect transversely at finitely many points.
The following are equivalent.
\begin{itemize}
\item 
$\alpha$ is isotopic relative to $\partial D \cup P$
to a simple closed curve which intersects $\beta$ at fewer points,
\item 
$\alpha$ and $\beta$ cobound a ``digon'',
that is, an embedded disc in $D \setminus P$
whose boundary consists of one subarc of $\alpha$
and one subarc of $\beta$.
\end{itemize}
\end{lem}

A proof can be found in
\cite[Proposition 3.7]{paris-rolfsen},
or \cite[Proposition 3.10]{flp}.

Suppose $\sigma \in \cH(D,P)$ is a homomorphism representing
an element of the kernel of the Lawrence-Krammer representation.
We will show that $\sigma$ is isotopic relative to $\partial D \cup P$
to the identity map.

Take $D$ to be the unit disc centered at the origin in the complex plane,
and take $p_1,\dots,p_n$ to be real numbers satisfying
$-1 < p_1 < \dots < p_n < 1$.
For $i = 1,\dots,n-1$,
let $E_i$ be the horizontal edge from $p_i$ to $p_{i+1}$.
For $i = 1,\dots,n$,
let $N_i$ be the noodle which winds around $p_i$ and no other punctures,
intersecting the real axis twice.
\begin{figure}
\centering
\begin{picture}(120,130)(0,-10)

\thinlines

\put(60,60){\ellipse{120}{120}}

\qbezier(50,1)(74,73)(83,70)
\qbezier(70,1)(92,67)(83,70)

\put(45,-10){$d_1$}
\put(65,-10){$d_2$}

\put(20,60){\circle*{4}}
\put(40,60){\circle*{4}}
\put(60,60){\circle*{4}}
\put(80,60){\circle*{4}}
\put(100,60){\circle*{4}}

\drawline(20,60)(40,60)
\end{picture}
\caption{The edge $E_1$ and the noodle $N_4$}
\label{stdnoodle}
\end{figure}
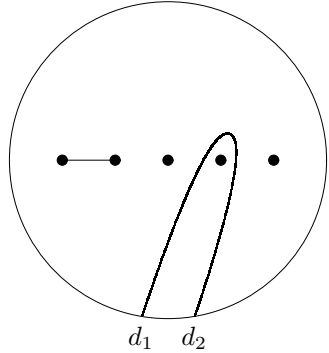
See Figure \ref{stdnoodle}.

Let $F$ be a fork such that $T(F) = E_1$.
Then $\langle N_3,F \rangle = 0$.
By the Basic Lemma, $\langle N_3,\sigma(F) \rangle = 0$.
By the Key Lemma, it follows that 
$\sigma(E_1)$ is isotopic relative to $\partial D \cup P$
to an arc which is disjoint from $N_3$.
By composing $\sigma$ with an element of $\cI(D,P)$ if necessary,
we can assume that $\sigma(E_1)$ is disjoint from $N_3$.

Similarly, $\sigma(E_1)$ can be isotoped so as to be disjoint from $N_4$.
By Lemma \ref{digon},
this isotopy can be performed by a sequence of moves
which consist of eliminating digons,
and hence do not introduce any new intersections with $N_3$.
Thus we can assume that 
$\sigma(E_1)$ is disjoint from both $N_3$ and $N_4$.

Continuing in this way,
we can assume that $\sigma(E_1)$ is disjoint from 
$N_i$ for all $i = 3,\dots,n$.
By applying one final isotopy relative to $\partial D \cup P$,
we can assume that $\sigma(E_1) = E_1$.
Note that we have not yet eliminated the possibility
that $\sigma$ reverses the orientation of $E_1$.

We can repeat the above procedure to isotope $\sigma(E_2)$ to $E_2$
while leaving $E_1$ fixed.
Continuing in this way,
we can assume that $\sigma(E_i) = E_i$
for all $i = 1,\dots,n-1$.
It follows that $\sigma$
must be isotopic relative to $\partial D \cup P$
to $(\Delta^2)^k$ for some $k \in \Z$,
where $\Delta^2$ is a Dehn twist about
the boundary of a collar neighborhood of $\partial D$.
However, the induced action of $\Delta^2$ on $H_2(\tilde{C})$
is simply multiplication by $q^{2n}t^2$.
Since $\sigma$ acts trivially on $H_1(\tilde{C})$
we have that $k=0$, so $\sigma$ represents the trivial braid.
\section{Matrices for the Lawrence-Krammer representation}
\label{lkmatrices}

In this section we give
an explicit description of the Lawrence-Krammer representation
in terms of matrices.

\begin{thm}
\label{lkmatrix}
$H_2(\tilde{C})$ is a free $\Lambda$-module of rank 
$\binom{n}{2}$.
There is a basis 
$$\{v_{j,k}:1 \le j < k \le n\}$$
on which the braid $\sigma_i$ acts as follows.
\[
\sigma_i(v_{j,k}) =
\left\{
\begin{array}{ll}
v_{j,k}
 &i \not \in \{j-1,j,k-1,k\}, \\
qv_{i,k}+(q^2-q)v_{i,j}+(1-q)v_{j,k}
 &i=j-1, \\
v_{j+1,k}
 &i=j \neq k-1, \\
qv_{j,i}+(1-q)v_{j,k}+(q^2-q)tv_{i,k}
 &i=k-1\neq j, \\
v_{j,k+1}
 &i=k \\
-tq^2 v_{j,k}
 &i=j=k-1.
\end{array}
\right.
\]
\end{thm}

We prove this theorem by
constructing a two-complex which is homotopy equivalent to $C$.
Our methods require some geometric intuition (read: ``hand-waving''),
and some details are left to the reader.

For $j = 1,\dots,n$, let $\xi_j$ be
a path in $D$ based at $d_1$ and passing counterclockwise around $p_j$,
and let $x_j$ be the arc $\{\xi_j,d_2\}$ in $C$.
Let $\tau_1$ be an arc from $d_1$ to $d_2$
and $\tau_2$ an arc from $d_2$ to $d_1$
such that $\tau_1\tau_2$ is 
a simple closed curve enclosing no puncture points.
Let $y$ be the arc $\{\tau_1,\tau_2\}$ in $C$.
Let $\cG = \{x_1,\dots,x_n,y\}$.

For $1 \le j \le n$, let
$$r_{j,j} = [x_j, y x_j y].$$
For $1 \le j < k \le n$, let
$$r_{j,k} = [x_j, y x_k y^{-1}].$$
Let $\cR = \{r_{j,k} : 1 \le j \le k \le n\}$.
Let $K$ be the Cayley complex of the presentation $\langle \cG|\cR \rangle$.
In other words, $K$ has one vertex,
one edge for each $g \in \cG$,
and one face $f_r$ for each $r \in \cR$,
where $f_r$ is attached to the $1$-skeleton
according to the word $r$.
We will show that $C$ is homotopy equivalent to $K$.

Let $\bar{C}$ be the set of 
{\em ordered} pairs of distinct points in $D \setminus P$.
This is the double cover of $C$
whose fundamental group is
normally generated by $x_1,\dots,x_n$ and $y^2$.

Let $X_j = y x_j y^{-1}$.
Let $Y = y^2$.
Let $\bar{\cG} = \{x_1,\dots,x_n,X_1,\dots,X_n,Y\}$.
For $1 \le j \le n$, let
\begin{eqnarray*}
\bar{r}_{j,j} &=& [x_j, X_j Y], \\
\bar{r}'_{j,j} &=& [X_j, Y x_j].
\end{eqnarray*}
For $1 \le j < k \le n$, let
\begin{eqnarray*}
\bar{r}_{j,k} &=& [x_j,X_k],\\
\bar{r}'_{j,k} &=& [X_j, Yx_kY^{-1}].
\end{eqnarray*}
Let 
\[
\bar{\cR} = 
\{\bar{r}_{j,k} : 1 \le j \le k \le n\} \cup
\{\bar{r}'_{j,k} : 1 \le j \le k \le n\}.
\]
Let $\bar{K}$ be 
the Cayley complex of $\langle \bar{\cG}|\bar{\cR} \rangle$.
Then $\bar{K}$ is homotopy equivalent to the double cover of $K$
whose fundamental group is
normally generated by $x_1,\dots,x_n$ and $y^2$.
To show that $C$ is homotopy equivalent to $K$,
it suffices to show that $\bar{C}$ is homotopy equivalent to $\bar{K}$.

Let $\pi \co \bar{C} \to D \setminus P$
be the map obtained by projection onto the first coordinate.
When restricted to the interior of $\bar{C}$,
this is a fiber bundle over the interior of $D \setminus P$
whose fiber is an $(n+1)$-times punctured open disc.

The base $D \setminus P$ is homotopy equivalent to a graph
with one vertex and $n$ edges
corresponding to $x_1,\dots,x_n$.
The fiber is homotopy equivalent to a graph
with one vertex and $n+1$ edges
corresponding to $X_1,\dots,X_n$ and $Y$.
The fiber bundle structure of $\bar{C}$
implies that it is homotopy equivalent to
the Cayley complex of a presentation
$\langle \bar{\cG}|\bar{\cR}' \rangle$,
where $\bar{\cR}'$ is a set of relations
equating $Y^{x_k}$ and $X_j^{x_k}$
to words in $\{X_1,\dots,X_n,Y\}$, for $j,k \in \{1,\dots,n\}$.
One can compute these relations $\bar{\cR}'$
by explicitly manipulating arcs in $\bar{C}$.
They are as follows.
\begin{eqnarray*}
Y^{x_k} &=& X_k Y X_k^{-1},\\
X_j^{x_k} &=&
  \left\{
  \begin{array}{ll}
    X_j Y X_j Y^{-1} X_j^{-1},&j=k\\
    X_k Y X_k^{-1} Y^{-1} X_j Y X_k Y^{-1} X_k^{-1},&j<k\\
    X_j & j>k. 
  \end{array}
  \right.
\end{eqnarray*}
One can transform the relations $\bar{\cR}'$ to $\bar{\cR}$
using moves which can be realized
by isotopy of the attaching maps of the faces
in the Cayley complex.
Thus $\bar{C}$ is homotopy equivalent to $\bar{K}$,
and hence $C$ is homotopy equivalent to $K$.

We are now ready to compute $H_2(\tilde{C})$.
Let $\cC_1$ and $\cC_2$ be the free $\Lambda$-modules
with bases $\{[g] : g \in \cG\}$ and $\{f_r:r \in \cR\}$ respectively.
For any word $w$ in $\cG$
we define $[w] \in \cC_1$ inductively 
according to the following rules
\begin{eqnarray*}
{[}1] &=& 0,\\
{[}gw] &=& [g] + \phi(g)[w],\\
{[}g^{-1}w] &=& \phi(g)^{-1}([w] - e_g),
\end{eqnarray*}
for any $g \in \cG$.
Then $H_2(\tilde{C})$
is the kernel of the map 
$\partial \co \cC_2 \to \cC_1$
given by $\partial f_r = [r]$.
We calculate the following.
\[
\partial f_r =
\left\{
\begin{array}{ll}
(1+q^{-1}t^{-1})((1-t)[x_j] + (q-1)[y]) &
\,\mbox{if}\,r = r_{j,j},\\ 
(q^{-1}-q^{-2})(-[x_j] + t[x_k] - (q-1)[y]) &
\mbox{ if }r = r_{j,k},\,\mbox{where}\,j<k.
\end{array}
\right.
\]
It is now an exercise in linear algebra to compute the kernel of this map.
It is a free $\Lambda$-module with bases
$\{v_{j,k} : 1 \le j < k \le n\}$, where
$$v_{j,k} = (q-1)f_{j,j} - (q-1)tf_{k,k} + (1-t)(1+qt)f_{j,k}.$$

We now define certain forks $F_{j,k}$
which will correspond to the basis vectors $v_{j,k}$.
Let $D$ be the unit disc centered at the origin in the complex plain.
Let $p_1,\dots,p_n$ lie on the real axis and satisfy
$-1<p_1 < \dots < p_n<1$.
Let $d_1$ and $d_2$ lie in the lower half plane,
with $d_1$ to the left of $d_2$.
For each $1 \le j < k \le n$,
let $F_{j,k}$ be a fork
which lies entirely in the closed lower half plane
such that the endpoints of $T(F)$ are $p_j$ and $p_k$.
Such an $F_{j,k}$ is uniquely determined up to isotopy by $j$ and $k$,
and will be called a {\em standard fork}.

Let $D' \subset D$ be a disc containing $F_{j,k}$
such that $D' \cap P = \{p_j,p_k\}$.
Let $C'$ be the set of unordered pairs of distinct points in $D'$.
Let $\tilde{C}'$ be the pre-image of $C'$ in $\tilde{C}$.
We can consider $v_{j,k}$ as an element of $H_2(\tilde{C}')$,
in which case it generates $H_2(\tilde{C}')$ as a $\Lambda$-module.
The surface $\tilde{\Sigma}_2(F_{j,k})$
lies in $\tilde{C}'$,
so must represent the homology class $\lambda v_{j,k}$
for some $\lambda \in \Lambda$.
The value of $\lambda$ does not depend on $j$ and $k$.
(Actually $\lambda = 1$, but we will not need this fact.)

To write $\sigma_i(v_{j,k})$ in terms of basis vectors,
we must find a $\Lambda$-linear combination of standard forks
which represents the same element of $H_2(\tilde{C})$
as the fork $\sigma_i(F_{j,k})$.

In the cases $i \not \in \{j-1,j,k-1,k\}$,
$i=j \neq k$,
and $i=k$,
there is no problem because
$\sigma_i(F_{j,k})$ is a standard fork.

In the case $i=j=k-1$,
the fork $\sigma_i(F_{j,k})$
has the same tine edge as $F_{j,k}$.
It follows that it represents the same surface in $\tilde{C}$,
up to a change in orientation
and application of a covering transformation.
With some thought,
or by pairing with an appropriate noodle,
it is not hard to check that
the correct formula is $\sigma_i(v_{j,k}) = -tq^2v_{j,k}$.

The remaining cases are
$i = j-1$ and $i=k-1 \neq j$.
We will use the following claim.

\begin{clm}
$\sigma_i(v_{j,k})$ is a linear combination of
basis vectors $v_{j',k'}$ which satisfy 
$j',k' \in \{i,i+1,j,k\}$.
\end{clm}

\begin{proof}
There exists a disc $D' \subset D$
such that $D'$ contains $\sigma(F_{j,k})$,
$D'$ contains $F_{j',k'}$
for all $j',k' \in \{i,i+1,j,k\}$ with $j' < k'$,
and $D' \cap P = \{p_i,p_{i+1},p_j,p_k\}$.
Let $C'$ be the set of
unordered pairs of distinct points in $D'$.
Let $\tilde{C}'$ be the pre-image of $C'$ in $\tilde{C}$.
Then $H_2(\tilde{C}')$
is a free $\Lambda$-module
with basis consisting of all $v_{j',k'}$
with $j',k' \in \{i,i+1,j,k\}$ and $j'<k'$.
But $\sigma(v_{j,k})$ can be considered as
an element of $H_2(\tilde{C}')$,
so must be a linear combination of these basis vectors.
\end{proof}

In the case $i=j-1$,
this claim implies that
$\sigma_i(F_{j,k})$ represents the same element of $H_2(\tilde{C})$
as some $\Lambda$-linear combination of the three standard forks
$F_{i,j}$, $F_{i,k}$, and $F_{j,k}$.
By pairing with some appropriate noodles
it is not hard to check that
the correct linear combination is as stated in Theorem \ref{lkmatrix}.
Similar methods can be used to
verify Theorem \ref{lkmatrix}
in the last remaining case, $i=k-1 \neq j$.
This completes the proof of Theorem \ref{lkmatrix}.

We conclude with some remarks on the
BMW representation of braid groups,
defined independently by Birman and Wenzl in \cite{birman-wenzl},
and by Murakami in \cite{murakami}.
V. Jones noticed a striking resemblance between
the matrices described in Theorem \ref{lkmatrix} 
and those of a certain irreducible summand of the BMW representation.
He asserted that the two representations should be the same
after some renormalization.
The details are worked out by Zinno in \cite{zinno}.
At present, there seems to be no deep explanation for this coincidence.
\providecommand{\bysame}{\leavevmode\hbox to3em{\hrulefill}\thinspace}


\begin{thebibliography}{FLP91}

\bibitem[Big99]{bigelow}
Stephen Bigelow, \emph{The burau representation is not faithful for $n=5$},
  Geometry and Topology \textbf{3} (1999), 397--404.

\bibitem[Bir74]{birman}
Joan~S. Birman, \emph{Braids, links, and mapping class groups}, Princeton
  University Press, Princeton, N.J., 1974, Annals of Mathematics Studies, No.
  82.

\bibitem[BW89]{birman-wenzl}
Joan~S. Birman and Hans Wenzl, \emph{Braids, link polynomials and a new
  algebra}, Trans. Amer. Math. Soc. \textbf{313} (1989), no.~1, 249--273.

\bibitem[FLP91]{flp}
A.~Fathi, F.~Laudenbach, and V.~Po\'enaru, \emph{Travaux de {T}hurston sur les
  surfaces}, Soci\'et\'e Math\'ematique de France, Montrouge, 1991, S\'eminaire
  Orsay, Reprint of {\it Travaux de Thurston sur les surfaces}, Soc.\ Math.\
  France, Paris, 1979 [MR 82m:57003], Ast\'erisque No. 66-67 (1991).

\bibitem[Kra]{krammer}
Daan Krammer, \emph{The braid group ${B}_4$ is linear}, (Preprint).

\bibitem[Law90]{lawrence}
R.~J. Lawrence, \emph{Homological representations of the {H}ecke algebra},
  Comm. Math. Phys. \textbf{135} (1990), no.~1, 141--191.

\bibitem[LP93]{long-paton}
D.~D. Long and M.~Paton, \emph{The {B}urau representation is not faithful for
  $n\geq 6$}, Topology \textbf{32} (1993), no.~2, 439--447.

\bibitem[Moo91]{moody}
John~Atwell Moody, \emph{The {B}urau representation of the braid group ${B}\sb
  n$ is unfaithful for large $n$}, Bull. Amer. Math. Soc. (N.S.) \textbf{25}
  (1991), no.~2, 379--384.

\bibitem[Mur87]{murakami}
Jun Murakami, \emph{The {K}auffman polynomial of links and representation
  theory}, Osaka J. Math. \textbf{24} (1987), no.~4, 745--758.

\bibitem[PR99]{paris-rolfsen}
L.~Paris and D.~Rolfsen, \emph{Geometric subgroups of surface braid groups},
  Ann. Inst. Fourier (Grenoble) \textbf{49} (1999), no.~2, 417--472.

\bibitem[Zin]{zinno}
Matthew~G. Zinno, \emph{{On Krammer's Representation of the Braid Group}},
  (arxiv: math.RT/ 0002136).

\end{thebibliography}
\end{document}